\documentclass[12 pt, twoside]{amsart}
\usepackage{amssymb}

\usepackage{graphicx}
\usepackage{amsmath,amssymb,latexsym}
\usepackage{amscd}

\newtheorem{thm}{Theorem}[section]

\newtheorem{lem}{Lemma}[section]
\begin{document}
\title[Non-existence of
unbounded Fatou components]{Non-existence of unbounded Fatou
components  of a meromorphic function}
\author[Zheng]{Zheng Jian-Hua$^{1}$}
\address{Department of Mathematical Sciences, Tsinghua University, Beijing, 100084,
P. R. China} \email{jzheng@math.tsinghua.edu.cn}

\author[Niamsup]{Piyapong Niamsup}
\address{Department of Mathematics, Faculty of Science, Chiang Mai University, Chiang
Mai, 50200, THAILAND} \email{scipnmsp@chiangmai.ac.th}
\thanks{$^{1}$Corresponding author} 
\maketitle


\begin{center}
    {\bf Abstract}
\end{center}

This paper is devoted to establish sufficient conditions under
which a transcendental meromorphic function has no unbounded Fatou
components and to extend some results for entire functions to
meromorphic function. Actually, we shall mainly discuss
non-existence of unbounded wandering domains of a meromorphic
function. The case for a composition of finitely many meromorphic
function with at least one of them being transcendental can be
also investigated in the argument of this paper.

\noindent {\bf Keywords and Phases.} Fatou set, Julia set

\section{Introduction and Main Results}

\indent Let $\mathcal{M}$ be the family of all functions
meromorphic in the complex plane $\mathbb{C}$ possibly outside at
most countable set, for example, a composition of finitely many
transcendental meromorphic functions is in $\mathcal{M}$. Here we
mean a function meromprphic in $\mathbb{C}$ with only one
essential singular point at $\infty$ by a transcendental
meromorphic function. We shall study iterations of element in
$\mathcal{M}$.

We denote the $n$th iteration of $f(z)\in \mathcal{M}$ by
$f^n(z)=f(f^{n-1}(z)),n=1,2,\ldots .$ Then $f^n(z)$ is well
defined for all $z \in \mathbb{C}$ outside a (possible) countable
set
$$E(f^n)=\bigcup_{j=0}^{n-1}f^{-j}(E(f)),$$
here $E(f)$ is the set of all essential
singular points of $f(z)$. Define the Fatou set $F(f)$ of $f(z)$
as
$$F(f)=\{z \in \bar{\mathbb{C}}: \{f^n(z)\}\ \textrm{is\ well\
defined}$$
$$\textrm{and\ normal\ in\ a\ neighborhood\ of}\ z\}$$ and
$J(f)=\bar{\mathbb{C}}\setminus F(f)$ is the Julia set of $f(z)$.
$F(f)$ is open and $J(f)$ is closed, non-empty and perfect. It is
well-known that both $F(f)$ and $J(f)$ are completely invariant
under $f(z)$, that is, $z \in F(f)$ if and only if $f(z) \in
F(f).$ And $F(f^n)=F(f)$ and $J(f^n)=J(f)$ for any positive
integer $n$. We shall consider components of the Fatou set $F(f)$
and hence let $U$ be a connected component of $F(f)$. Since $F(f)$
is completely invariant under $f$, $f^n(U)$ is contained in $F(f)$
and connected, so there exists a Fatou component $U_n$ such that
$f^n(U)\subseteq U_n$. If for some $n\geq 1, f^n(U)\subseteq U$,
that is, $U_n=U$, then $U$ is called a periodic component of
$F(f)$ and such the smallest integer $n$ is the period of periodic
component $U$. In particular, a periodic component of period one
is also called invariant. If for some $n$, $U_n$ is periodic, but
$U$ is not periodic, then $U$ is called pre-periodic; A periodic
component $U$ of period $p$ can be of the following five types:
(i) attracting domain when $U$ contains a point $a$ such that
$f^p(a)=a$ and $|(f^p)'(a)|<1$ and $f^{np}|_U\rightarrow a$ as
$n\rightarrow\infty$; (ii) parabolic domain when there exists a
point $a\in\partial U$ such that $f^p(a)=a$ and $(f^p)'(a)=e^{2\pi
i\alpha}$ for $\alpha\in\mathbb{Q}$ and $f^{np}|_U\rightarrow a$
as $n\rightarrow\infty$; (iii) Baker domain when $f^{np}|_U
\rightarrow a \in
\partial U \cup \{\infty\}$ as $n\rightarrow \infty$ and
$f^p(z)$ is not defined at $z=a$; (iv) Siegel disk when $U$ is
simply connected and contains a point $a$ such that $f^p(a)=a$ and
$\phi\circ f^p\circ\phi^{-1}(z)=e^{2\pi i\alpha}z$ for some real
irrational number $\alpha$ and a conformal mapping $\phi$ of $U$
onto the unit disk with $\phi(a)=0$; (v)Herman ring when $U$ is
doubly connected and $\phi\circ f^p\circ\phi^{-1}(z)=e^{2\pi
i\alpha}z$ for some real irrational number $\alpha$ and a
conformal mapping $\phi$ of $U$ onto $\{1<|z|<r\}$. $U$ is called
wandering if it is neither periodic nor preperiodic, that is,
$U_n\cap U_m =\emptyset$ for all $n\neq m.$ For the basic
knowledge of dynamics of a meromorphic function, the reader is
referred to \cite{Be} and the book \cite{Zheng}.

If for a function $f\in\mathcal{M}$, $f^{-2}(E(f))$ contains at
least three distinct points, then
$$J(f)=\overline{\bigcup\limits_{n=1}^\infty f^{-n}(E(f))},$$
and in any case, what we should mention is that for every $n\geq
1$, $f^n(z)$ is analytic on $F(f)$. In particular, this result
holds for a composition of finitely many meromorphic functions.

Our study in this paper relies on the Nevanlinna theory of value distribution.
To the end, let us recall some basic concepts and notations in the theory. Let $f(z)$ be a meromorphic function in $\mathbb{C}$. Define
$$m(r,f)=\int_0^{2\pi}\log^+|f(re^{i\theta})|d\theta$$
and
$$N(r,f)=\int_0^r{n(t,f)-n(0,f)\over t}dt+n(0,f)\log r,$$
where $n(t,f)$ is the number of poles of $f(z)$ in the disk
$\{|z|\leq t\}$, and
$$T(r,f)=m(r,f)+N(r,f)$$
which is known as the Nevanlinna characteristic function of
$f(z)$. The quantity $\delta (\infty, f)$ is the Nevanlinna
deficiency of $f$ at $\infty$, defined by the following formula
$$\delta (\infty, f)=\liminf_{r\rightarrow\infty}{m(r,f)\over T(r,f)}=1-\limsup_{r\rightarrow\infty}{N(r,f)\over T(r,f)}.$$
(See \cite{Hayman1}). The growth order and lower order of $f(z)$ are
defined respectively by
\begin{equation*}
    \lambda(f) = \limsup_{r\rightarrow+\infty} \frac{\log \ T(r,f)}{\log \ r}
\end{equation*}
and
\begin{equation*}
    \mu(f) = \liminf_{r\rightarrow+\infty} \frac{\log \ T(r,f)}{\log \ r}.
\end{equation*}

In this paper, we take into account the question, raised by I. N.
Baker in 1984, of whether every component of $F(f)$ of a
transcendental entire function $f(z)$ is bounded if its growth is
sufficiently small. Baker \cite{Baker1} shown by an example that
the order $1/2$ and minimal type is the best possible growth
condition in terms of order. Following I. N. Baker's question, a
number of papers gave some sufficient conditions which confirm
Baker's question for the case of entire functions.

Zheng \cite{Zheng2} made a discussion of non-existence of
unbounded Fatou components of a meromorphic function and actually
the method in \cite{Zheng2} is available in proving the following

\begin{thm}\label{thm1.1}
    Let $f(z)$ be a function in $\mathcal{M}$. If we
    have
    \begin{equation}\label{1.1}
        \limsup_{r\rightarrow +\infty}\frac{L(r,f)}{r} = +\infty,
    \end{equation}
where $L(r,f)= {\rm min} \{|f(z)|:|z|=r\}$, then the Fatou set,
$F(f)$, of $f$ has no unbounded preperiodic or periodic
components.
\\ In particular, $f$ has no Baker domains.
\end{thm}

Theorem \ref{thm1.1} confirms that an entire function whose growth
does not exceed order $1/2$ and minimal type has no unbounded
preperiodic or periodic components, whereas the result for the
case of order less than $1/2$ was proved in several papers, see
\cite{Stallard} and \cite{Hi1}. In view of a well-known result
that (\ref{1.1}) is satisfied for a transcendental meromorphic
function with lower order $\mu(f)<1/2$ and
$\delta(\infty,f)>1-\cos(\mu(f)\pi)$,
 Theorem \ref{thm1.1} also confirms that such a meromorphic function
has no unbounded preperiodic or periodic components. And it is
described by an example in Zheng \cite{Zheng2} that the condition
(\ref{1.1}) is sharpen. For a composition $g(z)=f_m \circ f_{m-1}
\circ \cdots \circ f_1(z)$ of finitely many transcendental
meromorphic functions $f_j(z)(j=1,2,\ldots,m;m\geq 1)$, from the
method of \cite{Zheng2} it follows that $F(g)$ has no unbounded
periodic or preperiodic components if for each $j$, there exits a
sequence of positive real numbers tending to infinity at which
$L(r,f_j)>r$ and (\ref{1.1}) holds for at least one $f_{j_0}$.

Therefore, the crucial point solving I. N. Baker's question is in
discussion of non-existence of unbounded wandering domains of a
meromorphic function. There are a series of results for the case
of entire functions on which some assumption on order less than
$1/2$ and the certain regularity of the growth are imposed. Let
$f(z)$ be an transcendental entire function with order $<1/2$.
Then every component of $F(f)$ is bounded, provided that one of
the following statements holds:

(1) ${\log M(2r,f)\over \log M(r,f)}\rightarrow c\geq 1\ \ {\rm
as}\ \ r\rightarrow\infty,$ (Stallard \cite{Stallard1}, 1993);

(2) ${\varphi'(x)\over\varphi(x)}\geq \frac{c}{x},$
for all sufficiently large $x$, where $\varphi(x)=\log M(e^x,f)$ and $c>1$
 (Anderson and Hinkkanen \cite{Hi1}, 1998);

(3) $\log M(r^m,f)\geq m^2\log M(r,f)$ for each $m>1$ and all sufficiently large $r$
(Hua and Yang \cite{HY}, 1999);

(4) $\mu(f)>0$ (Wang \cite{Wang}, 2001).

A straightforward calculation deduces that an entire function
satisfying the Stallard assumption with $c>1$ must be of lower
order at least $\log c/\log 2.$ However, an entire function with
$0<\mu\leq\lambda(f)<\infty$ must satisfy the Hua and Yang's
assumption for $m$ with $\mu(f)m>\lambda(f)$. In fact, choosing
$\varepsilon>0$ with $(\mu-\varepsilon)m>\lambda+2\varepsilon$, we
have for sufficiently large $r>0$
\begin{eqnarray}\label{1.2}
\log M(r^m,f)>(r^m)^{\mu-\varepsilon}>r^\varepsilon r^{\lambda+\varepsilon}
\geq r^\varepsilon\log M(r,f).\end{eqnarray}
What we should mention is that by modify a little the proof given in \cite{HY},
Hua and Yang's assumption for sufficiently large $m$ instead of each $m>1$ suffices
to confirm their result to be true.

Zheng and Wang \cite{Zheng3} in 2004 proved the following

\begin{thm}\label{thm1.3}\ \ Let $f(z)$ be a transcendental entire function. If there exists
a $d>1$ such that for all sufficiently large $r>0$ we can find a $\tilde{r}\in [r,r^d]$
satisfying
\begin{equation}\label{1.3}
\log L(\tilde{r},f)\geq d\log M(r,f),\end{equation}
then every component of $F(f)$ is bounded.\end{thm}

In \cite{Zheng3} they also made a discussion of the case of
composition of a number of entire functions. In 2005, Hinkkanen
\cite{Hi2} also gave a weaker condition than (\ref{1.3}), that is,
the coefficient "$d$" before $\log M(r,f)$ is replaced by
$"d(1-(\log r)^{-\delta})"$ with $\delta>0$.



In this paper, in view of the Nevanlinna theory of a meromorphic
function, we consider the case of a meromorphic function and our
main result is the following.

\begin{thm}\label{thm1.5}
    Let $f(z)$ be a transcendental
    meromorphic function and such that for some $\alpha \in
    (0,1)$ and $D>d>1$ and all the sufficiently large $r$, there
    exists an $t\in (r,r^d)$ satisfying
    \begin{equation}\label{1old}
        \log L(t,f) > \alpha T(r,f), j=1,2,\cdots,m.
    \end{equation}
    and
    \begin{equation}\label{n1.5}T(r^d,f)\geq DT(r,f).\end{equation}
    Then $F(f)$ has no unbounded components.
\end{thm}

Actually, the assumption in Theorem \ref{thm1.5} is also a
sufficient condition of existence of buried points of the Julia
set of a meromorphic function with at least one pole and which is
not the form $f(z)=a +(z-a)^{-p}e^{g(z)}$. For such a meromorphic
function, $J(f)=\overline{\bigcup_{j=0}^\infty f^{-j}(\infty)}$
and from Theorem \ref{thm1.5} $\infty$ is a buried point of $f(z)$
and therefore so are all prepoles.

As a consequence of Theorem \ref{thm1.5}, we have the following

\begin{thm}\label{thm1.4}
    Let $f(z)$ be a transcendental
    meromorphic function with
    $$\delta(\infty,f)>1-\cos(\pi\lambda(f))$$
    and $\lambda(f)<1/2$ and $\mu(f)>0$.
Then $F(f)$ has no unbounded components.
\end{thm}

In particular, Wang's result can be deduced from Theorem
\label{thm1.4}.

\section{The Proof of Theorems}

To prove Theorems, we need some preliminary results. First
preliminary result will be established by using the hyperbolic
metric and it has independent significance. To the end, let us
recall some properties on the hyperbolic metric, see
(\cite{Ahlfors}, \cite{Beardon}), etc. An open set $W$ in
$\mathbb{C}$ is called hyperbolic if $\mathbb{C}\setminus W $
contains at least two points (note $\infty$ has been kicked out of
$W$). Let $U$ be a hyperbolic domains\ in $\mathbb{C}$.
$\lambda_{U}(z)$ is the density of the hyperbolic metric on $U$
and $\rho_U(z_1, z_2)$ stands for the hyperbolic distance between
$z_1$ and $z_2$ in $U$, i.e.
\begin{equation*}
    \rho_U(z_1, z_2)=\inf_{\gamma \in U}
    \int_{\gamma}\lambda_U(z)|dz|,
\end{equation*}
where $\gamma$ is a Jordan curve connecting $z_1$ and $z_2$ in
$U$. For a hyperbolic open set $W$, the hyperbolic density
$\lambda_W(z)$ of $W$ is the hyperbolic density for each component
of $W$. Then we convent that the hyperbolic distance between two
points which are in disjoint components equals to $\infty$ and the
hyperbolic distance of two points $a$ and $b$ in one component $U$
equals to $\rho_W(a,b)=\rho_U(a,b)$. For a fixed point $a\not\in
W$, introduce a domain constant
$$C_W(a)=\inf\{|z-a|\lambda_W(z):z\in W\}.$$

If $U$ is simply-connected and $d(z, \partial U)$ is a
euclidean distance between $z \in U$ and $\partial U$, then for
any $z \in U$,
\begin{equation}\label{estimate_density}
    \frac{1}{2d(z, \partial U)}\leq \lambda_U(z) \leq \frac{2}{d(z, \partial
    U)}.
\end{equation}
Let $f:U \rightarrow V$ be analytic, where both $U$ and $V$ are
hyperbolic domains. By the principle of hyperbolic metric, we have
\begin{equation}\label{Priciple_of_hyperbolic}
    \rho_V(f(z_1), f(z_2))\leq \rho_U(z_1, z_2),\ {\rm for}\ z_1,
    z_2 \in U.
\end{equation}
In particular, if $U\subset V$, then $\lambda_V(z)\leq
\lambda_U(z)$ for $z\in U$.

\begin{lem}\label{lem1}(cf. Zheng \cite{Zheng})\ \ Let
$U$ be a hyperbolic domain and $f(z)$ a function such that each $f^n(z)$ is
analytic in $U$ and $\bigcup_{n=0}^\infty f^n(U)\subset W$. If for some fixed point $a\not\in W$, $C_W(a)>0$ and $f^n|_U\rightarrow
\infty$, then for any compact subset $K$ of $U$ there exists a positive constant $M=M(K)$ such that
\begin{equation}\label{1.8}M^{-1}|f^n(z)|\leq |f^n(w)|\leq M|f^n(z)|\ \ {\rm for}\ z, w\in K.\end{equation}
\end{lem}

{\bf Proof.}\ \  Under the assumption of Lemma \ref{lem1}, we obtain
    \begin{equation}\label{rho_U_n}
        \rho_{W}(z) \geq \frac{C_W(a)}{|z-a|} \geq \frac{C_W(a)}{|z|+|a|}.
    \end{equation}
It follows that
\begin{eqnarray}\label{lamma_U_n}
  \rho_{f^n(U)}(f^n(z), f^n(w)) &\geq& \rho_{W}(f^n(z), f^n(w))\nonumber\\
  &\geq& C_W(a)\left|\int_{|f^n(z)|}^{|f^n(w)|}\frac{dr}{r+|a|}\right| \nonumber \\
   &=& C_W(a)\left|\log\frac{|f^n(z)|+|a|}{|f^n(w)|+|a|}\right|.
\end{eqnarray}

Set $A= {\rm max}\{\lambda_{U}(z,w): z,\ w\in K \}.$ Clearly $A
\in (0, +\infty)$. From (\ref{Priciple_of_hyperbolic}), we have
\begin{equation}\label{eq_0f_A}
    \rho_{f^n(U)}(f^n(z), f^n(w)) \leq \rho_{U}(z, w)\leq A .
\end{equation}
Therefore, combining (\ref{lamma_U_n}) and
(\ref{eq_0f_A}) gives
\begin{equation}\label{contradict}
    |f^n(z)|+ |a|\leq (|f^n(w)|+ |a|)e^{A/C_W(a)}.
\end{equation}
This immediately completes the proof of Lemma \ref{lem1}.

The following is Lemma of Zheng \cite{Zheng2}( also see Theorem
1.6.7 of \cite{Zheng}).

\begin{lem}\label{lemn1} Let $f:U\rightarrow U$ map the hyperbolic domain
$U\subset \mathbb{C}$  analytically without fixed points and
without isolated boundary points into itself. If
$f^n|_U\rightarrow\infty (n\rightarrow\infty),$ then for any
compact subset $K$ of $U$, we have (\ref{1.8}) for some
$M=M(K)>0$.
\end{lem}

The following is a consequence of Lemma \ref{lem1} and Lemma
\ref{lemn1}, which is of independent significance.

\begin{thm}\label{thm2.1}\ \ Let $f(z)$ be a function in $\mathcal{M}$.
If $F(f)$ contains an unbounded component, then for any compact
subset $K$ of $F(f)$ with $f^n|_K\rightarrow\infty$ as
$n\rightarrow\infty$, we have a positive constant $M=M(K)$ such
that (\ref{1.8}) holds.
\end{thm}

{\bf Proof.} \ \ Assume without any loss of generalities that $K$
is contained in a component $U$ of $F(f)$. If $J(f)$ has one
unbounded component, then we can find a subset $\Gamma$ of $J(f)$
such that $\mathbb{C}\setminus\Gamma$ is simply-connected. Then in
view of Lemma \ref{lem1} we shall get $M=M(K)$ such that
(\ref{1.8}) holds by noting that $\bigcup_{n=0}^\infty
f^n(U)\subset W=\mathbb{C}\setminus\Gamma.$

Now assume that $J(f)$ only has bounded components and thus $F(f)$
has only one unbounded component denoted by $V$. If
$\bigcup_{n=0}^\infty f^n(U)$ does not intersect $V$, then in view
of the fact that $V$ has only bounded boundary components we can
choose a path $\Gamma$ in $V$ tending to $\infty$ such that
$\bigcup_{n=0}^\infty f^n(U)\subset W=\mathbb{C}\setminus\Gamma.$
Thus as we did above, the result of Theorem \ref{thm2.1} follows.

Let us consider the case when
$U\subseteq\bigcup_{n=0}^{\infty}f^{-n}(V)$. If $V$ is preperioidc
or periodic, then an application of Lemma \ref{lemn1} yields the
desired result of Theorem \ref{thm2.1}; If $V$ is wandering, then
for some $m>1$, $\bigcup_{n=m}^\infty f^n(U)$ does not intersect
$V$ and therefore we can prove Theorem \ref{thm2.1} in this case.

The second preliminary result comes from the Poisson formula.

\begin{lem}\label{lem3}\ \ Let $f(z)$ be meromorphic on $\{|z|\leq 3R\}$. Then there exists
a $r\in (R,2R)$ such that on $|z|=r$, we have
\begin{equation}\label{1.9}
\log^+ |f(z)|\leq  KT(3R,f).\end{equation} where $K(\leq 24)$ is a
universal constant, that is, it is independent of $R, r$ and
$f$.\end{lem}

{\bf Proof.}\ \ Set $D=\{|z|\leq\frac{5}{2}R\}$. We denote by
$G_D(\zeta,z)$ the Green function of $D$, that is,
$$G_D(\zeta,z)=\log\left|{(2.5R)^2-\bar{z}\zeta\over
2.5R(\zeta-z)}\right|, \ z, \zeta\in D.$$ A simple calculation
implies that
$$G_D(\zeta,z)\leq \log\frac{5R}{|\zeta-z|}$$ and for
$\zeta=2.5Re^{i\theta}$ and $r=|z|\leq 2R$,
$$\frac{\partial}{\partial\vec{n}}G_D(\zeta,z)ds={\rm
Re}{2.5Re^{i\theta}+z\over 2.5Re^{i\theta}-z}d\theta\leq{2.5R+r\over
2.5R-r}d\theta\leq 9d\theta.$$ In view of the Poisson formula, we
have
\begin{eqnarray*}
\log|f(z)|&=&\frac{1}{2\pi}\int_{\partial
D}\log|f(\zeta)|\frac{\partial}{\partial\vec{n}}G_D(\zeta,z)ds\\
&-&\sum_{a_n\in
D}G_D(a_n,z)+\sum_{b_n\in D}G_D(b_n,z)\\
&\leq& 9m(2.5R,f)+\sum_{b_n\in
D}\log\frac{5R}{|b_n-z|},\end{eqnarray*} where $a_n$ is a zero and
$b_n$ a pole of $f(z)$ in $D$ counted according to their
multiplicities. According to the definition of $N(r,f)$, we have
\begin{eqnarray*}
n(2.5R,f)&\leq&\left(\log\frac{6}{5}\right)^{-1}\int_{2.5R}^{3R}{n(t,f)\over
t}dt\\
&\leq& 6N(3R,f).\end{eqnarray*}

From the Boutroux-Cartan Theorem it follows that
$$\prod_{n=1}^N|z-b_n|\geq\left(\frac{R}{2e}\right)^N,\ \
N=n(2.5R,f),$$ for all $z\in \mathbb{C}$ outside at most $N$ disks
$(\gamma)$ the total sum of whose diameters does not exceed $R/2$.
Therefore there exists a $r\in [R,2R]$ such that
$\{|z|=r\}\cap(\gamma)=\emptyset$ and then on the circle $|z|=r$, we
have
$$\log^+|f(z)|\leq 9m(2.5R,f)+N\log 10e<24T(3R,f).$$

Thus we complete the proof of Lemma \ref{lem3}.

\

\noindent {\bf Proof Of Theorem \ref{thm1.5}.} For $\alpha>0$,
there exists a natural number $k$ such that $D^{k-1}\alpha\geq 1$.
Set $h=d^{k}$. In view of (\ref{1old}) and (\ref{n1.5}), for all
$r\geq R_0$, we have a $t\in (r^{d^{k-1}},r^h)$ such that
\begin{eqnarray}\label{1}
\log L(t,f)&\geq &\alpha T(r^{d^{k-1}},f) \geq  \alpha D^{k-1}
T(r,f)\nonumber\\ &\geq &T(r,f) ,\ \ {\rm on}\ |z|=t.
\end{eqnarray}

From Lemma \ref{lem3}, we have \begin{equation}\label{2} \log
|f(z)|\leq KT(3r,f),\ {\rm for}\ |z|\leq 2r,\end{equation} where
$K$ is a positive constant independent of $f$ and $r$.

Take a positive integer $m$ such that $D^{(m-1)k-1}>Kd^{mk}=Kh^m$.
Suppose that $f$ has an unbounded Fatou component, say $U$. Assume
that $U$ intersects $|z|=R_0$, otherwise we magnify $R_0$. Take a
point $z_0$ in $U\cap\{|z|=R_0\}$. Draw a curve $\gamma\in U$ from
$z_0$ to $U\cap\{|z|=R_0^H\},\ H=h^{m}$ such that $\gamma\subset
\{|z|=R^H_0\}$ except the end point of $\gamma$.

Then there exists a $z_1\in \gamma\cap\{R_0\leq |z|\leq 2R_0\}$
such that $\log |f(z_1)|\leq KT(3R_0,f)$. And there exists a
$r_1\in (R_0^{h^{m-1}},R_0^H)$ such that
\begin{eqnarray}\label{n2}
\log L(r_1,f)&\geq&
T(R_0^{h^{m-1}},f)=T(R_0^{d^{(m-1)k}},f)\nonumber\\
&\geq & D^{(m-1)k-1}T(R_0^d,f)>Kh^mT(3R_0,f),
\end{eqnarray}
on $|z|=r_1.$ Set $R_1=\exp(KT(3R_0,f_1))$. Then
\begin{equation}\label{1}
        f_1(\gamma)\cap\{|z|<R_1\}\not=\emptyset\ {\rm and}\
        f_1(\gamma)\cap\{|z|>R_1^{H}\}\not=\emptyset.
    \end{equation}
By the same argument as above, we have a $z_2\in
f(\gamma)\cap\{R_1\leq |z|\leq 2R_1\}$ such that $\log
|f(z_2)|\leq KT(3R_1,f)$ and a $r_2\in (R_1^{h^{m-1}},R_1^{H})$
such that
$$\log L(r_2,f)\geq {h^{m}}KT(3R_1,f) ,\ \ {\rm on}\
|z|=r_2.$$ Set $R_2=\exp(KT(3R_1,f_2))$. Then since the circle
$\{|z|=r_2\}$ intersects $f_1(\gamma)$, we have
\begin{equation}\label{2}
        f^2(\gamma)\cap\{|z|<R_2\}\not=\emptyset\ {\rm and}\
        f^2(\gamma)\cap\{|z|>R_2^H\}\not=\emptyset.
    \end{equation}
Define $R_n=\exp(KT(3R_{n-1},f))$ inductively. Then for each $n>0$
we always have
$$f^n(\gamma)\cap\{|z|<R_n \}\neq \emptyset$$ and
$$f^n(\gamma)\cap\{|z|\geq R_n^{H} \}\neq \emptyset.$$
Thus there is two points $z_{n},\ w_{n} \in \gamma$ such that

\begin{equation}\label{relation_of_g}
    |f^n(z_{n})| > R_{n}^H >|f^n(w_{n})|^H.
\end{equation}
Combining (\ref{relation_of_g}) and Theorem \ref{thm2.1} gives
\begin{equation}\label{contradict}
    |f^n(w_{n})|^H < |f^n(z_{n})|\leq M|f^n(w_{n})|.
\end{equation}
This is impossible as $n\rightarrow \infty,$ because  $a$ and
$e^{2A}$ are constants but $H>1$ and $|f^n(z_{n})| \rightarrow
+\infty $ as $n \rightarrow +\infty.$

This completes the proof of Theorem \ref{thm1.5}.

\

To prove Theorem 1.4, we need the following result, which was proved
by Gol'dberg and Sokolovskaya \cite{GS}.

\begin{lem}\label{lem2}    Let $f(z)$ be a transcendental
    meromorphic function with $\delta(\infty,f)>1-\cos(\pi\lambda(f))$
    and $\lambda(f)<1/2$. Then
    $$\underline{\log{\rm dens}}E>0,$$
    where $E=\{r>0:\log L(r,f)>\alpha T(r,f)\}$ for some positive
    $\alpha$.\end{lem}

In fact Lemma \ref{lem2} asserts that for sufficiently large
$r>0$, we can find a $t\in [r,r^d]$ for some $d>1$ with
$$\log L(t,f)>\alpha T(r,f).$$

For a function $f(z)$ with $0<\mu(f)\leq\lambda(f)<+\infty$, we
easily see that
$$\lim\limits_{r\rightarrow\infty}{T(r^d,f)\over T(r,f)}=\infty$$
 for $d$ with $d\mu(f)>\lambda(f)$.

 Therefore Theorem 1.4 follows immediately from Theorem
 \ref{thm1.5}.


\section{Conclusion}

By means of a careful calculation, indeed we can prove the
following result: a transcendental meromorphic function $f(z)$ has
no unbounded components of its Fatou set if for some $1<d<D$ and
all sufficiently large $r$ there exists a $t\in [r,r^d]$ such that
$$\log L(t,f)>DT(r,f).$$

The argument of this paper is also available in establishing the
corresponding  results for a composition of finitely many
meromorphic functions at least one of which is transcendental.

\end{document}